\documentclass[12pt,psamsfonts]{amsart}

\usepackage{amsmath,amssymb}
\usepackage{graphicx}
\usepackage[dvips]{psfrag}

\newtheorem{definition}{Definition}
\newtheorem{theorem}{Theorem}

\title[ Motion under attraction and rotation]{
On the  motion  under focal   attraction in a rotating medium}
\author[ Jorge Sotomayor]{}

\thanks{The  author is a fellow of CNPq and has the partial support of
CNPq Grant 473747/2006-5.}


\begin{document}
 \maketitle

\centerline{\scshape   Jorge Sotomayor}
\medskip

{\footnotesize \centerline{ Instituto de Matem\'atica e Estat\'\i
stica, Universidade de S\~ao Paulo} \centerline{Rua do Mat\~ao
1010, Cidade Universit\'aria} \centerline{05.508-090, S\~ao Paulo,
SP, Brasil } \centerline{\email{ sotp@ime.usp.br}}}

\medskip

\bigskip

\begin{quote}{\normalfont\fontsize{8}{10}\selectfont
{\bfseries Abstract.} New results  are established here  on the
phase portraits and bifurcations of
 the  kinematic   model  in system (\ref{eq:01}),
 first presented
  by
   H.K. Wilson  in \cite{wils}, and  by him attributed  to L. Markus
   (unpublished). A new, self-sufficient, study  which
   extends  that of  \cite{wils} and allows an essential  conclusion for the
   applicability of the model
   is reported here.
\par}
\end{quote}

\section{Introduction}\label{sec:1}

Consider the following  family of planar differential equations
depending on  three
 real parameters $(\omega, v, R)$, with $\omega \geq 0,  \;  v > 0,  \; R >
 0$,

\begin{equation}
\label{eq:01}
\begin{array}{lclr}
\dot{x} & = & -\omega y +
 v\frac{R - x}{\sqrt{(R - x)^2 + y^2}}, & \\
\dot{y} & = & \,\,\, \, \omega x - v\frac{y}{\sqrt{(R - x)^2 +
y^2}}.
\end{array}
\end{equation}

The {\it solutions}, also called {\it orbits}, of system
(\ref{eq:01}) describe the motion of certain {\it entities} (such
as particles or micro-organisms) attracted by a {\it focus} $F =
(R,0)$  (such as a light source or a  magnetic pole) toward which
they move with {\it velocity} $ v$; the motion takes place in a
{\it medium}\, (such as a fluid) which rotates with {\it angular
velocity} $\omega$  around a fixed point,  located  at the origin
$O=(0,0)$.

  The focal point $F =
(R,0)$, where system (\ref{eq:01}) is undefined,  will be regarded
as a {\it singular point}, outside of which it is analytic. A
point  at which the components of the system vanish will be
referred  to as an {\it equilibrium
 point}.

According  to  Wilson \cite {wils}, p. 297, this is a model
suggested by L. Markus (see \cite{wils}, p. viii,) for the motion
of {\it phototropic platyhelminthes}
---light-seeking flatworms--- swimming in a liquid filling a
shallow circular recipient with section
$$G_R = \{ x^2 + y^2 \leq R
\}. $$
 No printed
 reference  source  containing  Markus'   suggestion is given in \cite{wils}.

\begin {definition} Denote by $W^s (Q)$ the {\it basin of attraction}
of an equilibrium  or singular point
 $Q$ of (\ref{eq:01}). This is the set of points $p_0$ such that $\varphi (t, p_0 ) \to Q$
 as $t \to m_+ (p_0),$
the right extreme of the maximal interval of $\varphi (t, p_0 ) $.
The  {\it basin of repulsion}, $W^u (Q)$,  is defined analogously
for  $t \to m_- (p_0)$, referring to the left extreme of the
maximal interval.
\end {definition}

Notice that for  $Q=F$ the approach to $F$ by orbits of
(\ref{eq:01}) happens in finite time, that is,  $ m_+ (p_0)$ or  $
m_-(p_0)$ is finite.

In  \cite {wils}, p. 298,  using the  Poincar\'e - Bendixson
Theorem and the  Bendixson Negative  Criterium,  it  is proved the
existence and uniqueness of an {\it equilibrium point}  $P=
P_{(\omega,
  v, R)}$ attracting all positively
{\it complete}   semi-orbits $\varphi (t, p_0 )$  of system
({\ref{eq:01}})  starting at  points $ p_0 \in G_R \setminus F$.
Complete means  that $ \varphi (t, p_0 )$
 is  defined for all $ t \in [0, \infty). $  Denote such set of points $p_0$  by $G_{R,
 +}$.

For the terms in ODE's not  defined here, besides \cite{wils}, the
reader can profit from  consulting Chicone \cite{car}, among other
more recent up to date books.

Under the assumption
 $\omega
> v/R$,  $P$  is
 \begin{equation} \label{eq:02}
 P_{(\omega, v, R)} = ((v/{\omega})^2 /R, \sqrt{R^2-(v/\omega)^2
}(v/\omega R)).
\end{equation}
\noindent This  corrects   a misprint  in \cite{wils}, p. 298.

 The   theorem
below  improves  the  results    on this subject outlined in
\cite{wils}. An elementary  proof will be given in section
\ref{sec:2}.

\begin {theorem} \label{th:01}
 For all $\omega \geq 0$, the region $ G_R \setminus F$
is positively  invariant, in fact the radial component of system
(\ref{eq:01}) is negative on the complement of the closed disk $C$
of center  $(R/2, 0)$ and radius $R/2$.

\noindent It holds that

\begin{itemize}

 \item[1.]  For  $ 0 \leq \omega  \leq v/R $,  $F$   is a global attractor:  \,
  $W^s (F) = \mathbb R ^2 \setminus F .$

\item[2.] For $\omega > v/R $ there is a unique hyperbolic
attracting  equilibrium $P$ located at (\ref{eq:02}) whose basin
of  attraction contains  $W^u (F)$, the basin of repulsion of $F$,
itself  a regular  analytic curve contained in $G_R \setminus F.$
 \item[3.] Also, $W^s (F)$, the basin of attraction  of $F$
is  a
 regular analytic
 curve
disjoint from  $G_R$.

\noindent In particular  it holds that for $\omega > v/R $,
$$G_{R, +} = \, G_R \cap W^s (P)\, =\, G_R \setminus
F.$$

\end{itemize}

\end{theorem}

\section{Proof of Theorem \ref{th:01}} \label{sec:2}

Performing the change of variables and parameter
 rescaling
\begin{equation}
\label{eq:03} x= R \bar x,\,  y=R \bar y, \,  t=\bar t R/v ,  \,
\omega = \bar \omega v/R
\end{equation}
and then removing the bars, obtain  system   (\ref{eq:01}) with
$R\, =\, v\, =\, 1 :$
 \begin{equation}
\label{eq:01.1} \dot{x}  =   -\omega y + \frac{1 - x}{\sqrt{(1 -
x)^2 + y^2}}, \,\,\,\,\dot{y}  =   \omega x - \frac{y}{\sqrt{(1 -
x)^2 + y^2}}.
\end{equation}
Writing this  equation  in   polar coordinates centered at $F$,
given by
\begin{equation} \label{eq:05} x= 1-r\cos\theta, \, \; y = \,r\sin\theta, \,
\end{equation}
 \noindent  and multiplying both components  by $r$,  which amounts to
 rescaling again the time, obtain

\begin{equation}
\label{eq:07}
 \dot{\theta}  =  -\omega(r-\cos\theta ), \; \;
\dot{r}  = \,\,\, \, r(\omega\sin\theta -1).
\end{equation}

 Derivation of $z = x^2 +
y^2$ (which, in polar coordinates,  is $z= 1+r^2 -2r\cos\theta$)
in the direction  of equation (\ref{eq:01.1}) (which,  in polar
coordinates, has the same direction field as (\ref{eq:07})), gives
$z^\prime = -2r(r-\cos\theta)$.

 Clearly   $z^\prime $ is negative above the
curve $r=\cos\theta$, which is  the polar expression for the
border of the disk $C$.

Calculation of the equilibrium of (\ref{eq:07}) on the plane $r
>0$ gives the point $P$  whose polar coordinates
$(\theta_P , r_P) $
   satisfy
$$\sin\theta_P=1/\omega,\, \, r_P =
\cos\theta_P =  \sqrt{1- (1/\omega)^2}.$$

The expression, $(x_P , y_P)$,  of $P$ in cartesian coordinates
follows:

 $$x_P=1- \sqrt{1-(1/\omega)^2} \sqrt{1- (1/\omega)^2} = (1/\omega)^2
 , \, \, y_P=\sqrt{1- (1/\omega)^2}/\omega,$$

\noindent which in the original coordinates, before the change in
(\ref{eq:03}), reproduces  the expression in  (\ref{eq:02}).

Calculation of the divergence, $\sigma$,  and Jacobian, $\delta$,
of  (\ref{eq:07})   gives

$$\sigma = -1,    \, \, \, \delta = -\omega^2 + \omega \sin\theta +
\omega^2 r \cos\theta+\omega^2 \cos^2 \theta .  $$

Evaluation at $P$, i.e. at $(\theta_P , r_P ) $, gives $\delta =
\omega^2 -1$. This  shows that  $P$ is an attracting  hyperbolic
equilibrium point of {\it node} or {\it focus} type.

Calculation of the discriminant, $\sigma ^2 -4\delta$, for values
of $\omega \geq 1$ gives that the transition of $P$ from  a node
to a focus occurs  at $\sqrt{5}/2$, that is at $\omega = \sqrt{5}
v /2R$, in the original parameters.

This proves  the first general assertion and item 1 in the
theorem.

For $\omega >0$, system (\ref{eq:07}) has two equilibria; one is a
hyperbolic saddle at  $S_- = (-\pi/2, 0)$ with unstable separatrix
along the $\theta$-axis and  attracting eigenvector parallel  to
$(1,(2\omega+1)/\omega)$.

The other equilibrium, located at $S_+ = (\pi/2, 0),$ is a
hyperbolic node for $\omega < 1$. For $\omega
> 1$, $S_+ = (\pi/2, 0)$  is a hyperbolic saddle with stable
separatrix along the $\theta$-axis and with repelling eigenvector,
for positive eigenvalue $\omega -1,$  parallel  to $(1,(-2\omega +
1)/\omega)$.

 For the equilibrium at
$(\pi/2, 0)$,  the transition at $\omega = 1$, for increasing
$\omega$,  gives a cubic {\it pitchfork}  node to saddle
bifurcation. For $\omega >  1$, bifurcate two attracting nodal
equilibrium points, one to $r > 0$ and the other to $r < 0 $,
which capture the unstable separatrices of the saddle.  Only the
first one, located at $P$, is seen in the original $(x,y)$-plane.

\begin{figure}[htb]
\psfrag{L}{Saddle $S_- $} \psfrag{G}{Border of $G$}
\psfrag{C}{Border of $C$} \psfrag{R}{Saddle $S_+ $}
\begin{center}
\includegraphics[height=6cm]{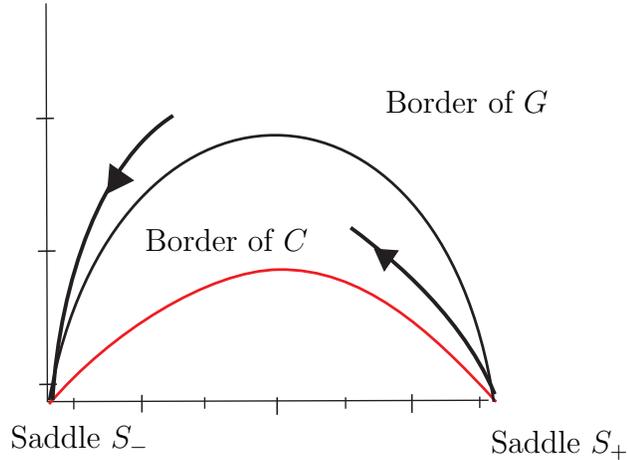}
\end{center}
\caption{Borders of $G$ and $C$ and Separatrix Directions
\label{fig1}}
\end{figure}

\begin{figure}[htb]
\psfrag{G}{Border of $G$} \psfrag{C}{Border of $C$}
\psfrag{F}{Focus $(R,0) $} \psfrag{P}{ $ P $}
\begin{center}
\includegraphics[height=6cm]{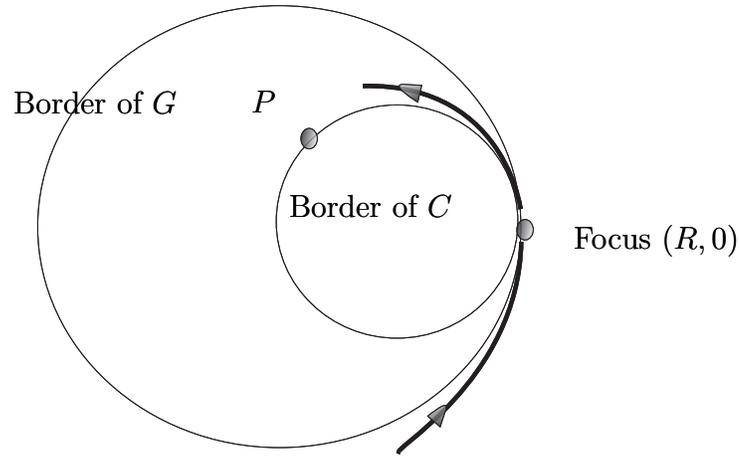}
\caption{Borders of $G$ and $C$, Equilibrium Point $P$ and
Separatrix Directions at Focus $F$ \label{fig2}}
\end{center}
\end{figure}

 Therefore,  the  slope  of the
unstable separatrix, at $S_+ $, is   $-2+ 1/\omega$  and that of
the stable separatrix at $S_- $ is $2+ 1/\omega$.

Comparison at $S_- $  of the slope the stable separatrix with the
slopes of  the borders of $G_R$ (for $R=1$), given by
$r=2\cos\theta$, which is $2$,  and of $C$, given by
$r=\cos\theta$,  which is $1$,  leads to the location  of $W^s
(F)$ outside of $G_R$.

The  location of $W^u (F)$ between $G_R$ and $C$ follows from
similar comparison of slopes, taking into account that at $S_+ $
the slopes of the borders of $G_R$ (for $R=1$) is $-2$,  and of
$C$  is $-1$.

 See Fig. \ref{fig1}, in polar coordinates, and Fig.
\ref{fig2}, in the original coordinates.

By continuity, for any $\omega > 1$ (that is greater than $v/R$, in
the  original coordinates),   the basin of repulsion of $F$ is
always contained in $W^s(P).$ This, for $\omega$ near $1$, is a
consequence of the structure of the local  {\it  saddle-node}
bifurcation. For large $\omega$  this follows from 
Poincar\'e-Bendixson Theorem and the  Bendixson Negative Criterium
\cite{wils}.

From the fact that $z^{\prime}$ is negative outside the disk $C$
follows that $W^s(P)$ coincides with $\mathbb R ^2 \setminus
W^s(F)$.

The two  main conclusions  in items 2 and 3  are established.

\section{Concluding Comments} \label{sec:3}

Item 3 is essential for the purpose of the biological  model in
\cite{wils}. In fact, from the knowledge of   the correct
location, at (\ref{eq:02}), of the equilibrium attracting all
entities moving with velocity $v$, their separation from similar
ones also contained in the positive invariant set $G_R$ but having
different characteristic velocities and, therefore, clustering at
different points of $G_R$, for large $t$. From this information
the removal of the entities for study in isolation  could   be
implemented practically. However,  no laboratory experiment report
where  the model has actually  been used is known by the author.

Both Wilson  \cite{wils} and  Sotomayor \cite{lic}, p. 264,
   the only written sources for system (\ref{eq:01}),
overlooked the property in item 3 of the theorem. Its
consideration was  proposed to the author toward 1995 by Dan Henry
(1944 - 2002).

\vskip 0.2cm

 \noindent{\bf Acknowledgement} The author is grateful to R.
Garcia for his help with the pictures and to L. F. Mello and A.
Gasull for their comments.

\end{document}